\documentclass[12pt]{article}

\usepackage{bm}
\usepackage{mathrsfs}
\usepackage{amsmath}
\usepackage{amsfonts}
\usepackage{mathtools}

\usepackage[margin=1in]{geometry}

\usepackage[english]{babel}
\usepackage{hyperref}
\usepackage{cleveref}

\usepackage{graphicx}
\usepackage{subcaption}
\usepackage{wrapfig}
\usepackage{float}
\usepackage{color}
\usepackage{verbatim}
\usepackage{algorithm}
\usepackage{algorithmic}
\usepackage{multirow}

\usepackage{wasysym}

\usepackage{amsopn}

\newcommand{\numberthis}{\addtocounter{equation}{1}\tag{\theequation}}
\newcommand{\del}{\partial}

\crefname{algorithm}{Algorithm}{Algorithm}

\crefname{section}{Section}{Sections}
\crefname{table}{Table}{Tables}
\crefname{figure}{Figure}{Figures}
\crefname{equation}{}{}
\crefrangeformat{equation}{\textup{(#3#1#4--#5#2#6)}}

\numberwithin{equation}{section}

\title{Density Propagation with Characteristics-based Deep Learning\thanks{This work has been submitted to IFAC for possible publication.}}

\author{Tenavi Nakamura-Zimmerer\thanks{Department of Applied Mathematics, Baskin School of Engineering, University of California, Santa Cruz (\texttt{tenakamu@ucsc.edu}, \texttt{venturi@ucsc.edu}, \texttt{qigong@soe.ucsc.edu}).}
\and
Daniele Venturi\footnotemark[2]
\and
Qi Gong\footnotemark[2]
\and
Wei Kang\thanks{Department of Applied Mathematics, Naval Postgraduate School, Monterey, CA (\texttt{wkang@nps.edu}).}}

\begin{document}

\maketitle

\begin{abstract}
Uncertainty propagation in nonlinear dynamic systems remains an outstanding problem in scientific computing and control. Numerous approaches have been developed, but are limited in their capability to tackle problems with more than a few uncertain variables or require large amounts of simulation data. In this paper, we propose a data-driven method for approximating joint probability density functions (PDFs) of nonlinear dynamic systems with initial condition and parameter uncertainty. Our approach leverages on the power of deep learning to deal with high-dimensional inputs, but we overcome the need for huge quantities of training data by encoding PDF evolution equations directly into the optimization problem. We demonstrate the potential of the proposed method by applying it to evaluate the robustness of a feedback controller for a six-dimensional rigid body with parameter uncertainty.
\end{abstract}

\section{Introduction}

We consider a $d$-dimensional system of nonlinear first-order ordinary differential equations,
\begin{equation}
\label{eq: system}
\dot {\bm x} = \bm f (\bm x) ,
\qquad
\bm x (0) = \bm x_0 (\xi) \sim \rho_0 (\bm x_0) ,
\end{equation}
where $\bm f (\bm x) : \mathbb R^d \to \mathbb R^d$ is the Lipschitz continuous vector field and $\rho_0 (\bm x_0) : \mathbb R^d \to [0, \infty)$ is a given finite-valued probability density function (PDF). This formulation includes uncertainty in parameter space if we augment the state with the set of uncertain parameters $\bm \mu$ and the vector field $\dot {\bm x} = \bm f (\bm x)$ with $\dot {\bm \mu} = \bm 0$. Randomness in the initial condition $\bm x_0 (\xi) $ induces randomness in the future state $\bm x (t)$. We can think of $\bm x (t)$ as a random process defined by the \emph{forward flow map}
\begin{equation}
\label{eq: flow map}
\bm x (t) = \bm \Phi (\bm x_0 (\xi), t) ,
\end{equation}
with a distribution characterized by the PDF $\rho (\bm x; t)$.

Since the vector field $\bm f (\bm x)$ is nonlinear, the flow map and thus the state PDF $\rho (\bm x; t)$ can be quite complicated. Often, statistical information such as the mean and variance of the state $\bm x (t)$ provide only poor characterizations of the distribution $\rho (\bm x; t)$ and the sensitivity of the state $\bm x (t)$ to perturbations in the initial conditions or parameters. Thus for many applications it is desirable to obtain a representation of the PDF itself.

Standard techniques for estimating $\rho (\bm x; t)$ include kernel density estimation (e.g. \cite{Scott2015}, Ch. 6), generalized polynomial chaos with flow map composition \cite{Luchtenburg2014}, probabilistic collocation \cite{Tatang1997}, PDF methods \cite{Cho2016}, and data-driven methods for conditional PDFs \cite{Brennan2018}. Many of these methods suffer from the curse of dimensionality, and even those that scale well to high dimensions are usually very data-hungry.

To fill this gap, we propose a computational method for approximating joint and conditional PDFs for systems where acquiring large amounts of sample data may be prohibitively expensive. Specifically, we model the time-dependent PDF with a neural network (NN), which we train by solving the PDF evolution equation (see \cref{sec: Liouville equation}) through a combination of least squares and method of characteristics approaches. This semi-supervised learning approach exploits both knowledge of the underlying dynamics and the ability to generate data. Once the NN is trained, it can predict the density at \emph{arbitrary points} in the computational domain at orders of magnitude more efficiently than numerical integration.

\section{The Liouville equation}
\label{sec: Liouville equation}

It is well-known (see e.g. \cite{Ehrendorfer1994}) that the PDF $\rho (\bm x; t)$ evolving according to the nonlinear dynamics \cref{eq: system} satisfies the \emph{Liouville transport equation},
\begin{equation}
\label{eq: Liouville equation}
\mathcal L [\rho] (\bm x, t) = 0 ,
\qquad
\rho (\bm x; 0) = \rho_0 (\bm x) ,
\end{equation}
for the linear operator $\mathcal L [ \cdot ]$ defined as
\begin{equation}
\label{eq: Liouville operator}
\mathcal L [\rho]
	\coloneqq \frac{\del}{\del t} \rho + \nabla_{\bm x} \cdot [\rho \bm f] .
\end{equation}
Equation \Cref{eq: Liouville equation} is a partial differential equation (PDE) in $d$ spatial dimensions plus time. Solving it directly requires a discretization of space and time which must be extremely fine because the support of the PDF can twist into thin curves within state space. Consequently, solving \cref{eq: Liouville equation} can be challenging even in low dimensions, and possibly intractable when the system is large.

On the other hand, using the method of characteristics one can derive a formal expression for the solution of \cref{eq: Liouville equation}:
\begin{equation}
\label{eq: Liouville solution}
\rho (\bm x; t) = \rho_0 (\bm x_0)
	\exp \left[ - \int_0^t \nabla_{\bm x} \cdot \bm f (\bm \Phi (\bm x_0, \tau)) d \tau \right] .
\end{equation}
Here $\bm x_0 = \bm \Phi_0 (\bm x, t)$, where $\bm \Phi_0 (\bm x, t)$ is the \emph{inverse flow map} which maps the state $\bm x (t)$ back to the initial condition $\bm x_0$ which generated it. For finite time $t$, this map is the unique inverse of the forward flow map \cref{eq: flow map} under standard smoothness assumptions on the vector field $\bm f (\bm x)$. This representation effectively decouples \emph{pointwise} solutions \cref{eq: Liouville equation}, allowing these to be computed independently. Such \emph{causality-free} characteristics-based methods have been used to solve various types of PDEs, such as scalar conservation laws \cite{Kang2017}, semilinear parabolic PDEs \cite{Han2018_PNAS}, and Hamilton-Jacobi-Bellman (HJB) equations \cite{Yegorov2018,Chow2019,Nakamura2019}.

To use \cref{eq: Liouville solution}, one must first find representations of the forward and inverse flow maps, which are high-dimensional nonlinear functions. Approximating these maps is perhaps more difficult than solving \cref{eq: Liouville equation} itself, and is a subject of ongoing research (see e.g. \cite{Luchtenburg2014, Lusch2018}). Still, \cref{eq: Liouville solution} allows us to evaluate the density along characteristics quite easily. In fact, every time a sample trajectory of \cref{eq: system} is computed, we can simultaneously propagate $\rho (\bm x; t)$ by
\begin{equation}
\label{eq: dp/dt}
\dot \rho = - \rho \, \nabla_{\bm x} \cdot \bm f (\bm x) ,
\end{equation}
which can obtained by rearranging terms in \crefrange{eq: Liouville equation}{eq: Liouville operator}. Having such density data readily available suggests a data-driven approach, but one which is augmented by knowledge of the underlying physics, i.e. that $\rho (\bm x;t)$ obeys the Liouville equation \cref{eq: Liouville equation}.

\section{Deep learning for probability density function approximation}
\label{sec: PINNs for PDF}

Deep learning offers an efficient way to approximate high-dimensional PDFs. Deep NNs are well-known for their ability to approximate arbitrary high-dimensional nonlinear functions given sufficient training data, and are the focus of much recent research in scientific computing. They are also computationally efficient, allowing the estimation of joint and conditional densities at millions of spatio-temporal coordinates in seconds. This in turn allows computation of statistics like means and covariances in high dimensions through Markov Chain Monte Carlo (MCMC) algorithms \cite{Brooks2011}.

In this paper, we model the desired PDF $\rho (\cdot)$ with fully-connected feedforward NNs, which we denote by $\rho_{\bm \theta} (\cdot)$ for a parameterization $\bm \theta$. Thus we have
$$
\rho (\bm x; t)
	\approx \rho_{\bm \theta} (\bm x; t)
	= g_L \circ g_{L-1} \circ \dots \circ g_1 (\bm x, t) ,
$$
where each layer $g_l (\cdot), l = 1, \dots L$, is defined as
$$
g_l (\bm y) = \sigma_l (\bm W_l \bm y + \bm b_l )
$$
for weight matrices $\bm W_l$, bias vectors $\bm b_l$, and nonlinear activation functions $\sigma_l (\cdot)$. The trainable parameters are $\bm \theta = \{ \bm W_l, \bm b_l \}_{l=1, \dots, L}$. In this paper we use $\sigma_l (\cdot) = \tanh (\cdot)$ for all hidden layers $l = 1, \dots, L-1$, and a linear output layer (i.e. $\sigma_L (\cdot)$ is the identity map).

In a na\"{i}ve application of deep learning, we are given a large set of training data $\mathcal D$ consisting of inputs $\{ \bm x^{(i)}, t^{(i)} \}_{i = 1, \dots, | \mathcal D |}$ and outputs $\{ \rho^{(i)} \}_{i = 1, \dots, | \mathcal D |}$, where $\rho^{(i)} \coloneqq \rho ( \bm x^{(i)}; t^{(i)} )$. Each data point corresponds to a point in the time series of a sample trajectory, of which there are many. To train the NN, we solve a nonlinear regression problem over $\bm \theta$ to match the NN predictions $\rho_{\bm \theta} ( \bm x^{(i)}; t^{(i)} )$ with the training data $\rho^{(i)}$.

Unfortunately, deep NNs notoriously require enormous quantities of data to train, but we consider the case where data is not necessarily abundant. This situation can arise when numerical integration of a system is expensive, or if we are building a model from real-world measurements of a system with known structure but uncertain parameters. To overcome a relative lack of data, we adapt the physically-motivated machine learning strategies proposed by \cite{Nakamura2019} for solving high-dimensional HJB equations in optimal control, and \cite{Raissi2019} and \cite{Sirignano2018} for other PDEs. We outline this process in the following sections.

We also note that \cite{Yang2019} recently proposed a generative-adversarial strategy to train probabilistic surrogate models of systems like \cref{eq: system}. This allows one to sample from $\rho_{\bm \theta} (\bm x; t)$, whereas our method predicts the density itself.

\subsection{Physics-informed learning}
\label{sec: physics-informed learning}

In \cite{Nakamura2019}, the authors improve data-efficiency in training by augmenting the regression loss with an additional term which encourages the NN to learn the gradient of the solution. In the context of optimal control, this gradient is computed naturally for each sample trajectory. No such simple relationship exists for the PDF evolution equation. Thus for physically-motivated regularization we turn to the least squares methods proposed by \cite{Raissi2019} and \cite{Sirignano2018}. Simply put, we would like $\rho_{\bm \theta} (\cdot)$ to satisfy the Liouville equation \cref{eq: Liouville equation}. To this end we define the PDE residual
\begin{equation}
\label{eq: Liouville residual}
R \left[ \rho_{\bm \theta} \right] (\bm x, t)
	\coloneqq \left( \mathcal L \left[ \rho_{\bm \theta} \right] (\bm x, t) \right)^2 ,
\end{equation}
which of course is everywhere zero if $\rho_{\bm \theta} \equiv \rho$. The partial derivatives $\del \rho_{\bm \theta} / \del t$ and $\nabla_{\bm x} \rho_{\bm \theta}$ appearing in \cref{eq: Liouville residual} can be calculated using automatic differentiation. This standard feature of machine learning software packages allows efficient computation of \emph{exact} partial derivatives, and is one of the essential drivers of the success of deep learning.

We minimize the residual over a set of collocation points, $\mathcal C \coloneqq \left \{ \bm x^{(j)}, t^{(j)} \right \}_{j = 1, \dots, |\mathcal C|}$. Such collocation points may include the training data as well as randomly sampled points which require no numerical integration to generate. Approximating the PDE at collocation points can alleviate the need for large amounts of training data, as well as generating a representation of the PDF which is guided by the underlying physics \cref{eq: Liouville equation}, rather than just regression. The main difference between our method and least squares approaches is that we exploit the ability to generate data along the characteristics of the PDF. We find that this two-pronged strategy is more effective than direct minimization of the PDE residual and boundary conditions.

We now introduce the \emph{physics-informed learning problem},
\begin{equation}
\label{eq: PINN loss}
\begin{array}{cc}
\underset{\bm \theta}{\text{minimize}} & L (\bm \theta; \mathcal D, \mathcal C)
	\coloneqq \text{loss}_{\rho} \left( \bm \theta; \mathcal D \right) + \lambda \cdot \text{loss}_{\mathcal L} \left( \bm \theta; \mathcal C \right) .
\end{array}
\end{equation}
Here $\text{loss}_{\mathcal L} \left( \bm \theta; \mathcal C \right)$ is the mean square PDE residual,
\begin{equation}
\label{eq: PDE loss}
\text{loss}_{\mathcal L} (\bm \theta; \mathcal C) \coloneqq
\frac{1}{|\mathcal C|} \sum_{j = 1}^{| \mathcal C |} R \left[ \rho_{\bm \theta} \right] \left( \bm x^{(j)}, t^{(j)} \right) .
\end{equation}
For the regression term, $\text{loss}_{\rho} (\bm \theta; \mathcal D)$, we find that a weighted mean square error on the log PDF is effective:
\begin{align*}
&\text{loss}_{\rho} (\bm \theta; \mathcal D) 
\numberthis
\label{eq: data loss}
 \coloneqq \frac{1}{|\mathcal D|} \sum_{i=1}^{| \mathcal D |} w_i \left[ \log \rho_{\bm \theta} \left( \bm x^{(i)}; t^{(i)} \right) - \log \rho^{(i)} \right]^2 ,
\end{align*}
for some weights $\{ w_i \}_{i=1, \dots, |\mathcal D|}$. Examples include
$$
w_i = \rho^{(i)} ,
\qquad
w_i = \sqrt{\rho^{(i)}} ,
\qquad
\text{and}
\qquad
w_i = 1 ,
$$
for $i = 1, \dots, |\mathcal D|$. The choice of appropriate weight depends on the scaling of the problem. In the implementation, we actually have the NN predict the log PDF; then to obtain $\rho_{\bm \theta} (\bm x; t)$ we exponentiate the NN output. Therefore the NN model naturally \emph{preserves positivity of the PDF}.

Lastly, $\lambda$ is a scalar weight which must be carefully chosen to balance the impact of the PDE residual term \cref{eq: PDE loss} with the regression loss \cref{eq: data loss}, which serves as a \emph{boundary condition} in the present context. Good choices of $\lambda$ are highly problem-dependent, and for some problems like the rigid body problem in \cref{sec: rigid body}, we have found success with increasing $\lambda$ over the course of training.

\subsection{Adaptive sample size selection}
\label{sec: adaptive learning}

The machine learning problem introduced in \cref{sec: physics-informed learning} is atypical in that data generation may be expensive for stiff, chaotic, or high-dimensional systems. Data scarcity can be alleviated by use of randomly sampled collocation points, but this is often not enough to obtain good results. In addition, randomly sampled collocation points might completely miss the support of the PDF, which can be twisted into small manifolds in the higher-dimensional space. For example, see the double gyre in \cite{Luchtenburg2014} or \cref{fig: kraichnan orszag pdf evolution} in this paper. It is worth mentioning that PDFs like this are poorly characterized by simple statistics such as means and variances.

On the other hand, we have the freedom to generate additional data in parallel to training, hence we can treat this as a semi-online learning problem, incorporating additional data as it becomes available. For simplicity, in the current implementation we train the NN in multiple rounds, allowing the optimizer to converge in each round. In between rounds we generate additional data as needed, which we use to continue training the model.

Standard early stopping tests in machine learning rely on generating extra data sets for measuring generalization accuracy, which can be problematic if data is scarce in the first place. Further, such tests do not provide estimates of how much data is needed for convergence. To address these issues, we modify the convergence test and sample size selection scheme developed by \cite{Nakamura2019} to the present problem.

Let $\mathcal D_r$ and $\mathcal C_r$ denote the training data and collocation sets available in the $r$th training round. The method proposed in \cite{Nakamura2019} works for a single data set, so we apply it to each loss term individually\footnote{We will postpone extending this idea to deal with loss functions with multiple data sets future work, but note that it can actually be advantageous to adjust each data set independently.}. Thus we consider $\mathcal D_r$ to be sufficiently large if
\begin{equation}
\label{eq: norm test D}
\frac{\sum_{m=1}^{| \bm \theta |} \text{Var}_{\mathcal D_r} \left[ \frac{\del \text{loss}_{\rho}}{\del \theta_m} \left( \bm \theta; \left( \bm x^{(i)}, t^{(i)} \right) \right) \right]}
	{\left| \mathcal D_r \right| \left \Vert \nabla_{\bm \theta} \text{loss}_{\rho} \left( \bm \theta; \mathcal D_r \right) \right \Vert_1}
	\leq \epsilon_{\rho} ,
\end{equation}
where $\epsilon_{\rho} > 0$ is a scalar parameter and $\text{Var}_{\mathcal D_r} [\cdot]$ denotes the \emph{sample} variance\footnote{Computing the variance for large data sets can become expensive, so when necessary we evaluate this term over a subset of the data.},
\begin{align*}
& {\textstyle \text{Var}_{\mathcal D_r} \left[ \frac{\del \text{loss}_{\rho}}{\del \theta_m} \left( \bm \theta; \left( \bm x^{(i)}, t^{(i)} \right) \right) \right]} \\
\coloneqq& {\textstyle \text{Var}_{\left( \bm x^{(i)}, t^{(i)} \right) \in \mathcal D_r} \left[ \frac{\del \text{loss}_{\rho}}{\del \theta_m} \left( \bm \theta; \left( \bm x^{(i)}, t^{(i)} \right) \right) \right] .}
\end{align*}
Similarly for $\mathcal C_r$ we require
\begin{equation}
\label{eq: norm test C}
\frac{\sum_{m=1}^{| \bm \theta |} \text{Var}_{\mathcal C_r} \left[ \frac{\del \text{loss}_{\mathcal L}}{\del \theta_m} \left( \bm \theta; \left( \bm x^{(i)}, t^{(i)} \right) \right) \right]}
	{\left| \mathcal C_r \right| \left \Vert \nabla_{\bm \theta} \text{loss}_{\mathcal L} \left( \bm \theta; \mathcal C_r \right) \right \Vert_1}
	\leq \epsilon_{\mathcal L} .
\end{equation}

Intuitively, the convergence tests \crefrange{eq: norm test D}{eq: norm test C} check if the sample gradients are good approximations of the ``true'' gradients taken over infinitely large data sets. Satisfying \crefrange{eq: norm test D}{eq: norm test C} do not indicate if the trained NN is a good model, only that feeding it more data would likely not improve it significantly. The main advantage of these tests is that if they are not satisfied, then they immediately provides estimates for \emph{how much} data should be generated for the next round, namely
\begin{equation}
\label{eq: sample size selection D}
s_{\rho} | \mathcal D_r | \geq \left| \mathcal D_{r+1} \right| \geq
	\frac{\sum_{m=1}^{| \bm \theta |} \text{Var}_{\mathcal D_r} \left[ \frac{\del \text{loss}_{\rho}}{\del \theta_m} \left( \bm \theta; \left( \bm x^{(i)}, , t^{(i)} \right) \right) \right]}
	{\epsilon_{\rho} \left \Vert \nabla_{\bm \theta} \text{loss}_{\rho} \left( \bm \theta; \mathcal D_r \right) \right \Vert_1} ,
\end{equation}
and similarly for $\mathcal C_{r+1}$,
\begin{equation}
\label{eq: sample size selection C}
s_{\mathcal L} | \mathcal C_r | \geq \left| \mathcal C_{r+1} \right| \geq
	\frac{\sum_{m=1}^{| \bm \theta |} \text{Var}_{\mathcal C_r} \left[ \frac{\del \text{loss}_{\mathcal L}}{\del \theta_m} \left( \bm \theta; \left( \bm x^{(i)}, , t^{(i)} \right) \right) \right]}
	{\epsilon_{\mathcal L} \left \Vert \nabla_{\bm \theta} \text{loss}_{\mathcal L} \left( \bm \theta; \mathcal C_r \right) \right \Vert_1} .
\end{equation}
In the above, $s_{\rho}, s_{\mathcal L} > 0$ are scalar parameters which prevent the sample size from growing too fast. For details we refer the reader to \cite{Nakamura2019}.

\subsection{Transfer learning for extending time horizons}
\label{sec: transfer learning}

For some problems where the density changes quickly over time or where a longer time horizon is desired, it can be difficult to directly learn the distribution over all $t \in [0, t_f]$. In such cases, we employ a simple \emph{transfer learning} approach to gradually extend the time horizon learned by the NN. To do this, we pick a monotone increasing sequence of time horizons $\{ t_k \}_{k = 1, \dots, N_t}$, where $0 < t_1 < t_f$ and $t_{N_t} = t_f$. We first train the NN to model the density over $t \in [0, t_1]$, then re-train the NN over $t \in [0, t_2]$. We continue in this way until the NN learns the full time horizon, $t \in [0, t_f]$.

\subsection{Model validation}

The final component of any machine learning algorithm is model validation. While finding theoretical convergence properties remains an active area of research (see e.g. \cite{Han2018_arXiv,Sirignano2018,Zhang2019}), our data-driven approach naturally allows for empirical testing of model accuracy. To do this, we first construct a \emph{validation} data set, $\mathcal D_{\text{val}}$, of trajectories generated from initial conditions drawn independently of those used for training. We evaluate each trajectory in $\mathcal D_{\text{val}}$ at the same sequence of time snapshots, $\{ t_k \}_{k=0, 1, \dots, K}$, where $t_0 = 0$ and $t_K = t_f$. Then for each $t_k$ we estimate the normalized root mean square error (NRMSE) of the model predictions:
\begin{equation}
\label{eq: validation}
\text{NRMSE} (\rho_{\bm \theta} ; t_k) \coloneqq \frac{\sqrt{\sum_{i = 1}^{| \mathcal D_{\text{val}} |} \left[ \rho_{\bm \theta} \left( \bm x^{(i)} ; t_k \right) - \rho \left( \bm x^{(i)} ; t_k \right) \right]^2}}{\sqrt{\sum_{i = 1}^{| \mathcal D_{\text{val}} |} \left[ \rho \left( \bm x^{(i)} ; t_k \right) \right]^2}} .
\end{equation}
Here $\bm x^{(i)} \coloneqq \bm x^{(i)} (t_k)$ is the $k$th snapshot along the $i$th trajectory. This simple validation framework provides a characterization of the NN's performance over the time horizon, which can be compared between target densities of different scale.

We summarize the full training procedure in \cref{alg: training}.

\begin{algorithm}[h!]
\caption{Multi-round characteristics-based learning}
\label{alg: training}
\begin{algorithmic}[1]
\STATE{Generate $\mathcal D_1$, $\mathcal D_{\text{val}}$ by sampling from $\rho_0 (\bm x_0)$ and integrating \cref{eq: system,eq: dp/dt} forward}
\FOR{time horizon $t_1, t_2, \dots, t_{N_t}$}
	\FOR{$r = 1, 2, \dots$}
		\STATE{Solve \cref{eq: PINN loss} to update $\bm \theta$}
		\IF{\cref{eq: norm test D,eq: norm test C} are satisfied}
			\STATE{\textbf{return} optimized $\bm \theta$, NRMSE \cref{eq: validation}}
		\ELSE
			\STATE{Integrate additional samples to satisfy \cref{eq: sample size selection D}}
			\STATE{Sample collocation points to satisfy \cref{eq: sample size selection C}}
		\ENDIF
	\ENDFOR
\ENDFOR
\end{algorithmic}
\end{algorithm}

\section{Results}

\subsection{Kraichnan-Orszag problem}
\label{sec: kraichnan-orszag}

First we study the Kraichnan-Orszag problem from \cite{Orszag1967,Wan2005}. The state dynamics are
\begin{equation}
\label{eq: kraichnan orszag dynamics}
\begin{array}{lcccr}
\dot x_1 = x_1 x_3 ,
&\quad&
\dot x_2 = - x_2 x_3 ,
&\quad&
\dot x_3 = - x_1^2 + x_2^2 .
\end{array}
\end{equation}
We consider independently normally distributed initial conditions such that the joint PDF $\rho_0 (\bm x_0)$ straddles the ``stochastic discontinuity'' on the $x_2 = 0$ axis:
\begin{equation}
\label{eq: kraichnan orszag initial condition}
\begin{dcases*}
x_1 (0) = x_{1,0} (\xi) \sim \mathcal N (1, 1/4^2) ,\\
x_2 (0) = x_{2,0} (\xi) \sim \mathcal N (0, 1/2^2) ,\\
x_3 (0) = x_{3,0} (\xi) \sim \mathcal N (0, 1/2^2) .
\end{dcases*}
\end{equation}
The Kraichnan-Orszag system is a canonical test problem in uncertainty quantification. We seek to model the time-evolving PDF for $t \in [0, 10]$, and compare the performance of our method using different optimization and collocation point selection strategies.

Using TensorFlow \cite{Tensorflow}, we implement a fully-connected feedforward NN with four hidden layers with 64 neurons each. For the weights in the regression loss term \cref{eq: data loss}, we pick $w_i = \rho^{(i)}$, $i = 1, \dots, | \mathcal D |$. We consider three strategies to train the network:
\begin{enumerate}
\item \emph{Adam:} We train for one round with a fixed data set, $\mathcal D$. We set $\lambda = 1/2$ and optimize using Adam \cite{Kingma2014} for a pre-fixed number of epochs and with manually tuned learning rates. We randomly sample a mini-batch of collocation points at each iteration. This implementation is similar to the ``Deep Galerkin Method'' (DGM) proposed by \cite{Sirignano2018}, but in the context of solving \cref{eq: Liouville equation} with trajectory data serving as a boundary condition.
\item \emph{Fixed L-BFGS:} We train for one round with a fixed data set, $\mathcal D$. We construct $\mathcal C$ as the union of $\mathcal D$ and an equally-sized set of uniformly sampled collocation points, fixed during training. We set $\lambda = 1/2$ and optimize using L-BFGS \cite{Byrd1995}. This implementation is similar to the method proposed by \cite{Raissi2019}.
\item \emph{Adaptive L-BFGS:} W use the multi-round model refinement scheme outlined in \cref{alg: training}, optimizing with L-BFGS.
\end{enumerate}

The first two implementations use a data set $\mathcal D$ of 500 sample trajectories evaluated at $K = 80$ time snapshots each, for a total of the $| \mathcal D | = 4 \cdot 10^4$ training data. The adaptive learning implementation starts from a smaller data set, $\mathcal D_1$, composed of 250 trajectories from the fixed data set, but we add new data each round. We choose $\mathcal C_1$ to be the union of $\mathcal D_1$ and an equally-sized set of uniformly sampled collocation points, also adding more points progressively. Lastly we set $\lambda = 1/2$, $s_{\rho} = s_{\mathcal L} = 2$, $\epsilon_{\rho} = 6 \cdot 10^{-4}$, and $\epsilon_{\mathcal L} = 3 \cdot 10^{-4}$. For all implementations we directly learn the density evolution over $t \in [0, 10]$.

\begin{figure}[t!]
\centering
\includegraphics[width = 0.6 \textwidth]{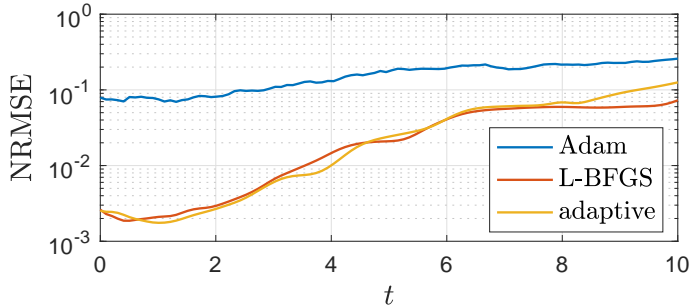}
\caption{Validation error with different training strategies when modeling a PDF evolving according to the Kraichnan-Orszag dynamics.}
\label{fig: ko error}
\end{figure}

To validate the models, we construct a second data set $\mathcal D_{\text{val}}$ of 500 sample trajectories evaluated at $K = 100$ time snapshots each, for a total of $| \mathcal D_{\text{val}} | = 5 \cdot 10^4$ validation data. \Cref{fig: ko error} shows the NRMSE \cref{eq: validation} evaluated on $\mathcal D_{\text{val}}$ for each of the implementations tested. From this plot, it is clear that training for a single round with L-BFGS is more effective than with Adam and collocation points randomly selected at each iteration. This was true for all hyperparameter configurations tested with Adam. In \cref{fig: ko error}, we see that the resulting NN is just as accurate as the other training strategies, even though it was trained without full access to a large data set. This suggests that the training strategies discussed in \cref{sec: PINNs for PDF} can be useful for solving data-scarce problems.

\begin{figure*}[t!]
\centering
\includegraphics[width = 0.328 \textwidth]{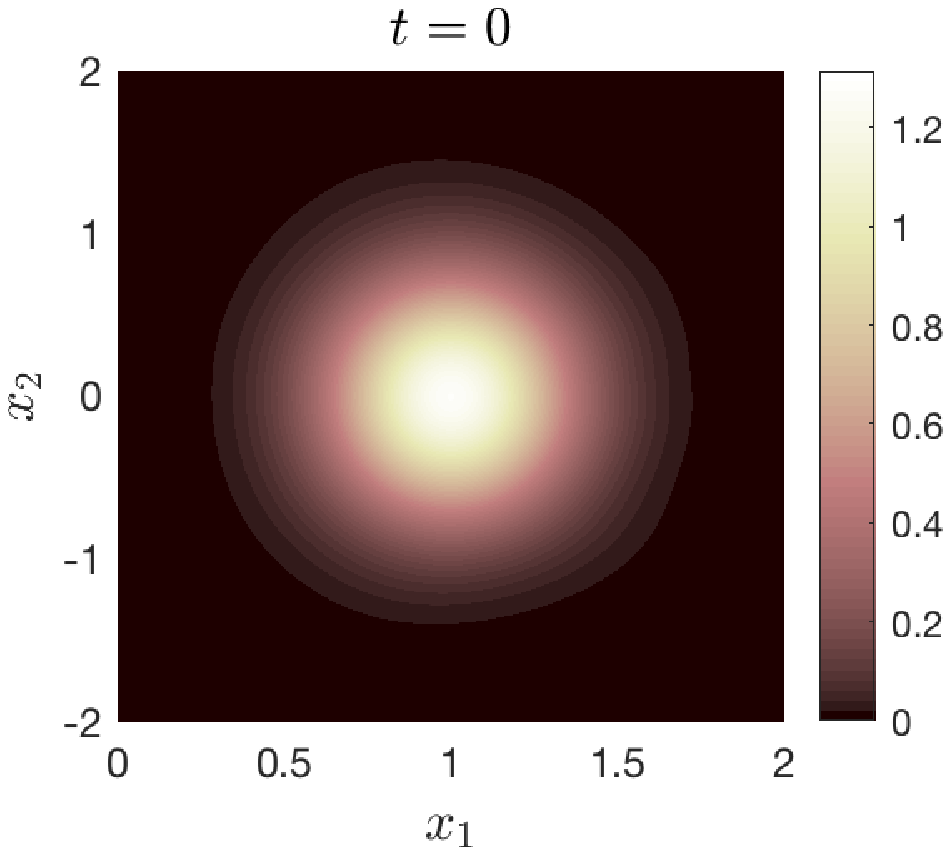}
\includegraphics[width = 0.328 \textwidth]{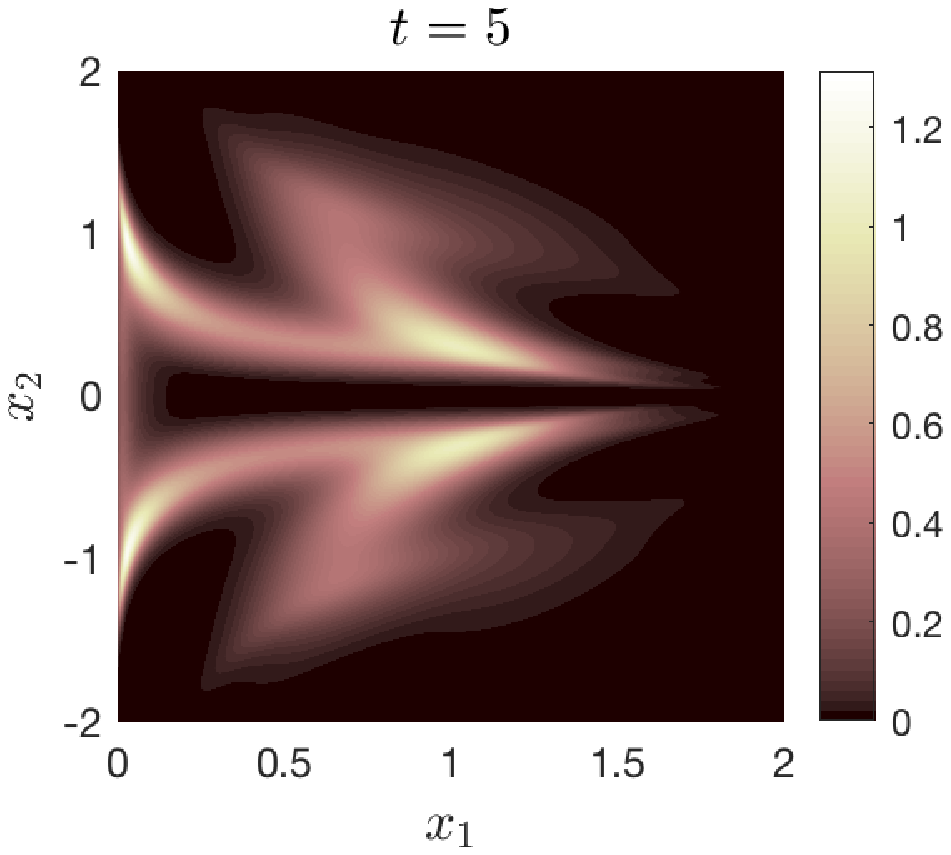}
\includegraphics[width = 0.328 \textwidth]{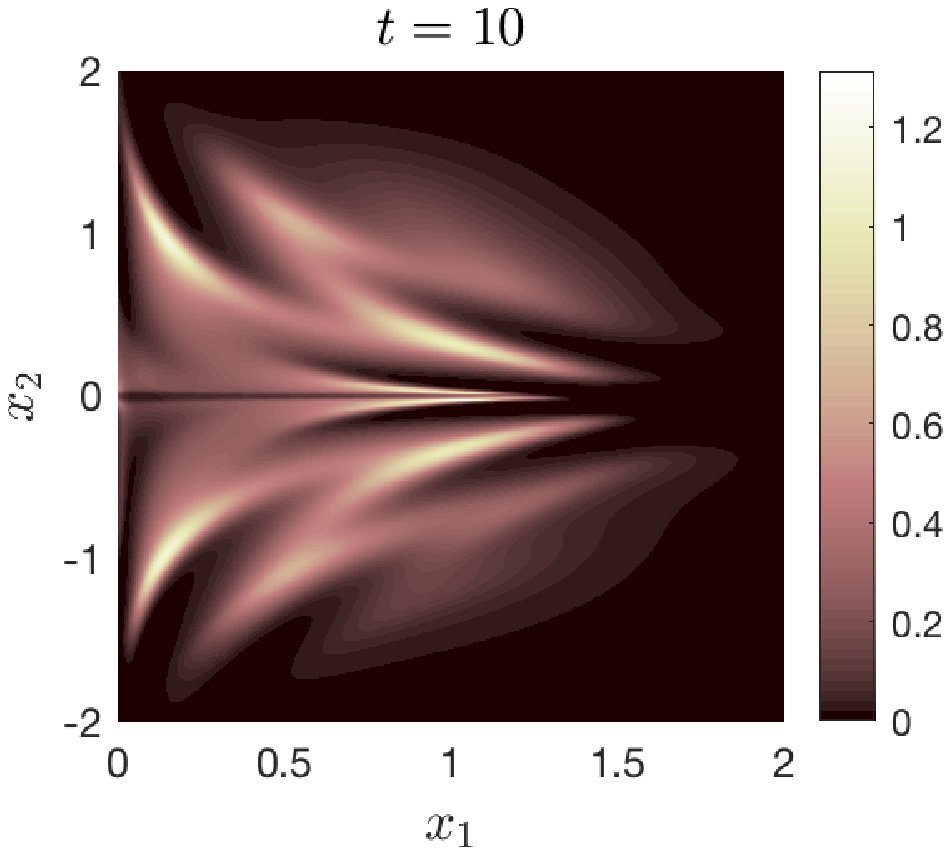}
\caption{Predicted marginal PDF, $\rho_{\bm \theta} (x_1, x_2; t)$, of the Kraichnan-Orszag system.}
\label{fig: kraichnan orszag pdf evolution}
\end{figure*}

Lastly in \cref{fig: kraichnan orszag pdf evolution} we plot the predicted marginal density, $\rho_{\bm \theta} (x_1, x_2; t) \approx \rho (x_1, x_2; t)$. The NN predicts the joint density on a tensor grid, which we then marginalize by quadrature integration. While dense grid integration is impractical in higher dimensions, here we use it only to visualize the complexity of the PDF captured by the NN.

\subsection{Application to quantifying controller robustness}
\label{sec: rigid body}

In this section we demonstrate the potential of the proposed method for use in evaluating controllers. We consider a six-dimensional rigid body model of a satellite controlled by reaction wheels. The state is written as $\bm x = (\bm v , \bm \omega)$, where $\bm v = (\phi, \theta, \psi)^T$ is the attitude of the satellite represented in Euler angles (roll, pitch, yaw) and $\bm \omega = (\omega_1, \omega_2, \omega_3)^T$ denotes the angular velocity in the body frame.  In this representation, the state dynamics are
\begin{equation}
\label{eq: satellite dynamics}
\begin{dcases*}
\dot{\bm v} = \bm E (\bm v) \bm \omega ,\\
\bm J \dot{\bm \omega} =  \bm S (\bm \omega) \bm R (\bm v) \bm h + \bm B \bm u ,
\end{dcases*}
\end{equation}
where $\bm E (\bm v), \bm S (\bm \omega), \bm R (\bm v) : \mathbb R^3 \to \mathbb R^{3 \times 3}$ are nonlinear matrix-valued functions\footnote{See \cite{Nakamura2019} for details.}; $\bm J \in \mathbb R^{3 \times 3}$ is the inertia matrix; $\bm h \in \mathbb R^3$ is the constant angular momentum of the system; and $\bm B \in \mathbb R^{3 \times 3}$ is a constant matrix. The system is controlled by applying a torque $\bm u (\bm v, \bm \omega) : \mathbb R^3 \times \mathbb R^3 \to \mathbb R^3$. We let
$$
\bm J = \begin{pmatrix}
2 & 0 & 0 \\
0 & 3 & 0 \\
0 & 0 & 4
\end{pmatrix} ,
\quad
\bm B = \begin{pmatrix}
\beta & 0.1 & 0.2 \\
0.2 & \beta & 0.3 \\
0.3 & 0.2 & \beta
\end{pmatrix} ,
\quad
\text{and}
\quad
\bm h = \begin{pmatrix}
1 \\
1 \\
1
\end{pmatrix} ,
$$
for a random parameter $\beta = \beta (\xi)$.

To stabilize the system, we design a linear quadratic regulator (LQR) by solving a linearized infinite horizon optimal control problem,
\begin{equation}
\label{eq: satellite OCP}
\left \{
\begin{array}{cl}
\underset{\bm u (\cdot)}{\text{minimize}} & \displaystyle \int_0^\infty W_1 \Vert \bm v \Vert^2 + W_2 \Vert \bm \omega \Vert^2 + W_3 \Vert \bm u \Vert^2 dt, \vspace{3pt} \\
\text{subject to} & \dot {\bm x} = \bm F \bm x + \bm G \bm u , \vspace{3pt}
\end{array}
\right .
\end{equation}
where $\bm F \in \mathbb R^{6 \times 6}$, $\bm G \in \mathbb R^{6 \times 3}$ are constant matrices obtained by linearizing \cref{eq: satellite dynamics} about $\bm x = (\bm v,  \bm \omega) = \bm 0$ and using a \emph{nominal} parameter value $\bar \beta = 1$. We choose the weights as $W_1 = 4$, $W_2 = 1/2$, and $W_3 = 8$. It is well-known that the minimizer of \cref{eq: satellite OCP} is a linear state feedback controller, $\bm u = \bm u_{\text{LQR}} (\bm x) = - \bm {K x}$, and that the gain matrix $\bm K$ is computed by solving a continuous algebraic Riccati equation (offline).

\begin{figure}[t!]
\centering
\includegraphics[width = 0.6 \textwidth]{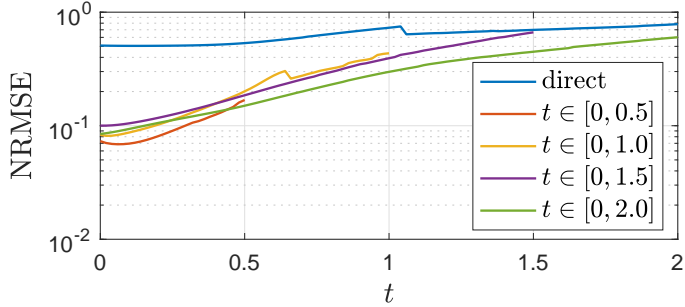}
\caption{Validation error with transfer learning and direct solution when modeling a PDF evolving according to the rigid body dynamics with LQR control.}
\label{fig: satellite error}
\end{figure}

Using the proposed density propagation framework, we model the evolution of a density of initial states over $t \in [0, 2]$, subject to the nonlinear dynamics \cref{eq: satellite dynamics} under LQR control. The system is augmented with an additional state $\beta$ for which $\dot \beta = 0$, bringing the total dimension to seven. Suppose that at $t=0$ the states are independently distributed according to
\begin{equation}
\label{eq: satellite initial conditions}
\begin{dcases*}
\phi_0 (\xi), \theta_0 (\xi), \psi_0 (\xi) \sim \mathcal N (0, (\pi/6)^2) , \\
\omega_{1,0} (\xi), \omega_{2,0} (\xi), \omega_{3,0} (\xi) \sim \mathcal N (0, 2^2) , \\
\beta (\xi) \sim {\textstyle \frac{1}{2}} \mathcal N (1/3, 1/9^2) + {\textstyle \frac{1}{2}} \mathcal N (1, 1/9^2) .
\end{dcases*}
\end{equation}
The bimodal Gaussian mixture distribution for $\beta$ is chosen to allow us to simultaneously study the behavior when $\beta$ is near the nominal value $\bar \beta$, and a scenario where $\beta < \bar \beta$, representing estimation error or actuator inefficiency.

\begin{figure*}[t!]
\centering
\begin{subfigure}{\textwidth}
\centering
\includegraphics[width = 0.328 \textwidth]{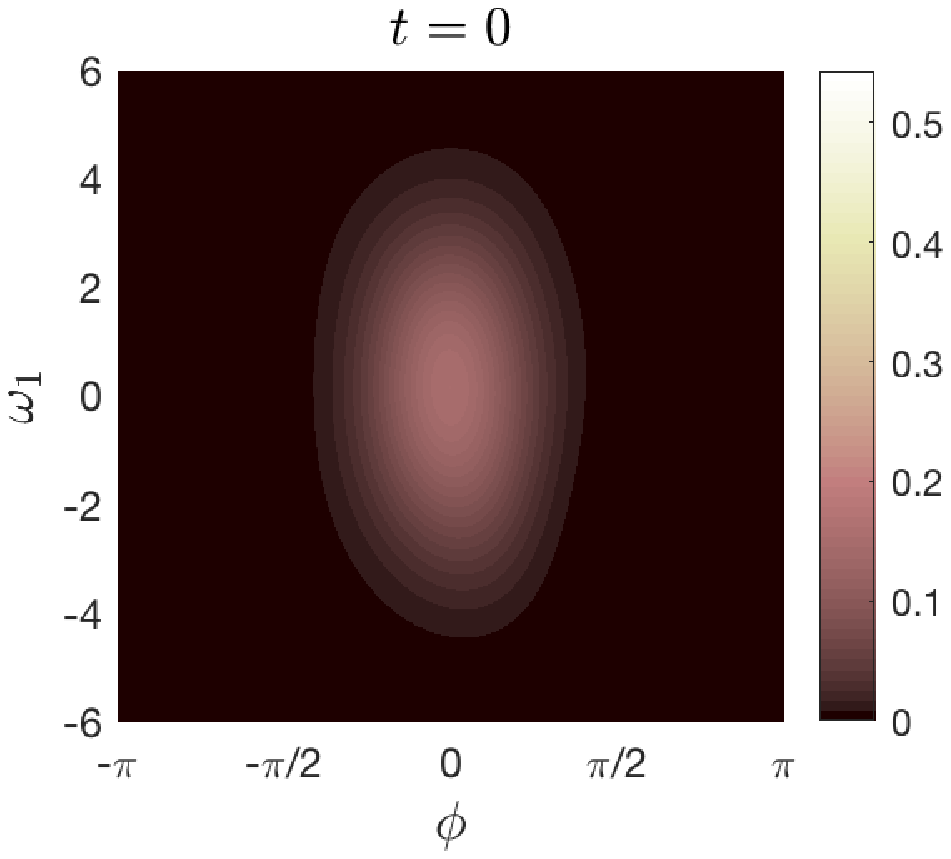}
\includegraphics[width = 0.328 \textwidth]{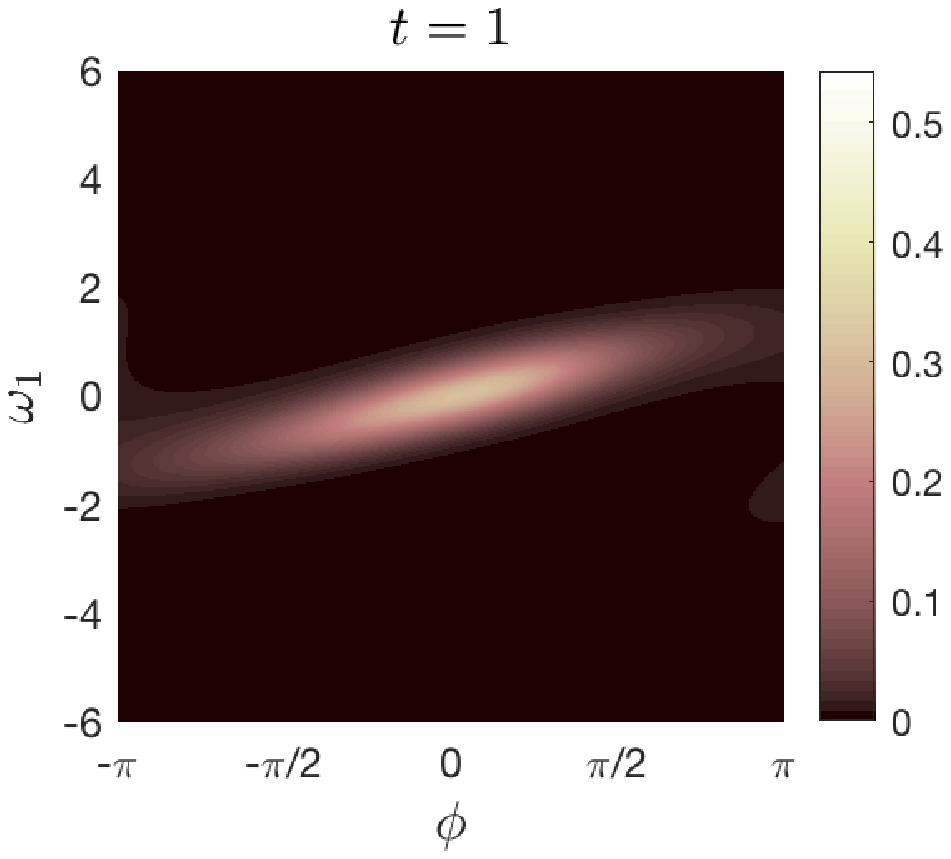}
\includegraphics[width = 0.328 \textwidth]{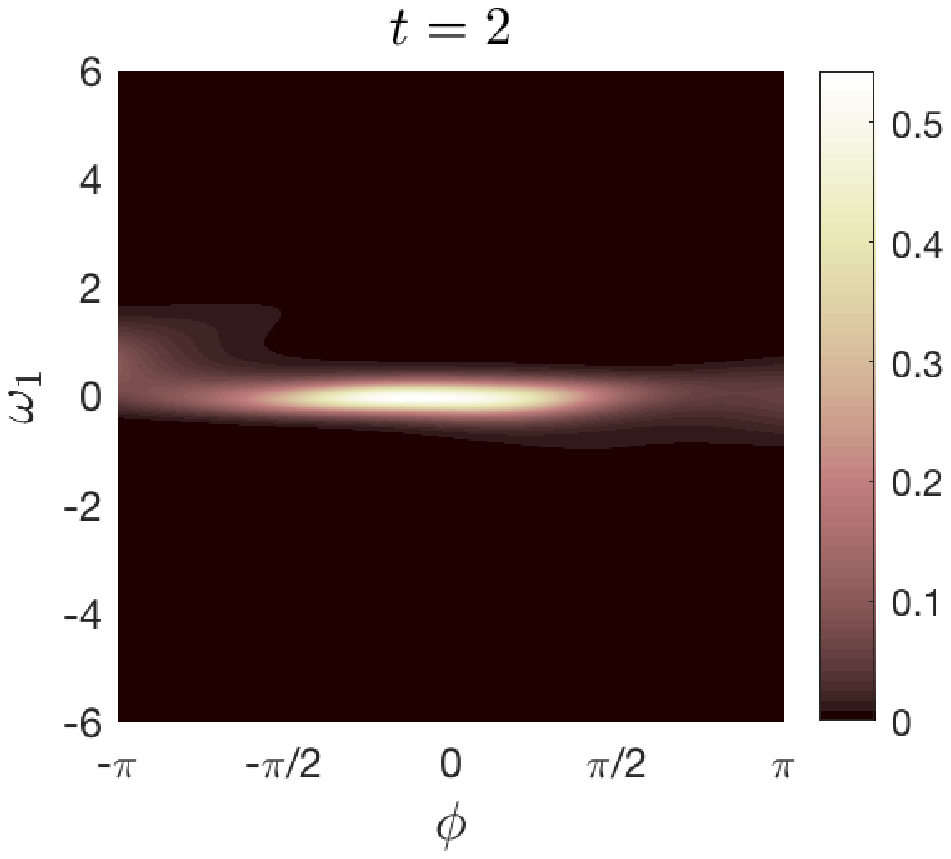}
\caption{$\rho_{\bm \theta} \left( \phi, \omega_1 \left| \theta = \psi = \omega_2 = \omega_3 = 0, \beta = \bar \beta = 1 \right. ; t\right)$}
\label{fig: satellite pdf 1 evolution}
\end{subfigure}
\begin{subfigure}{\textwidth}
\centering
\includegraphics[width = 0.328 \textwidth]{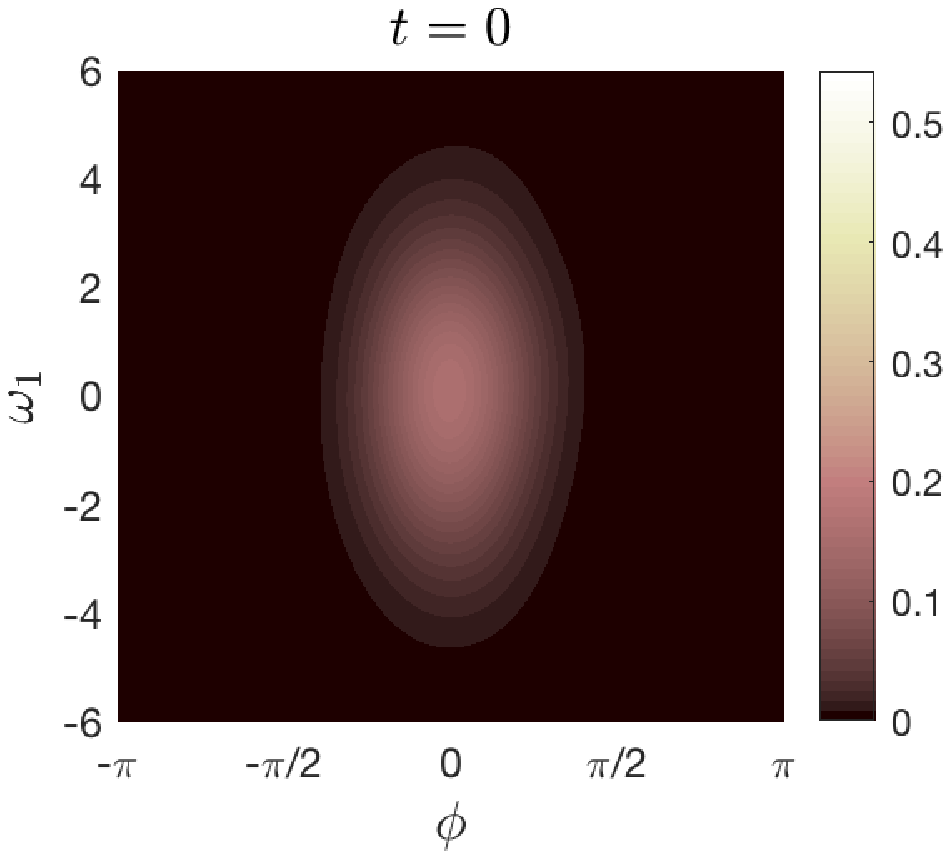}
\includegraphics[width = 0.328 \textwidth]{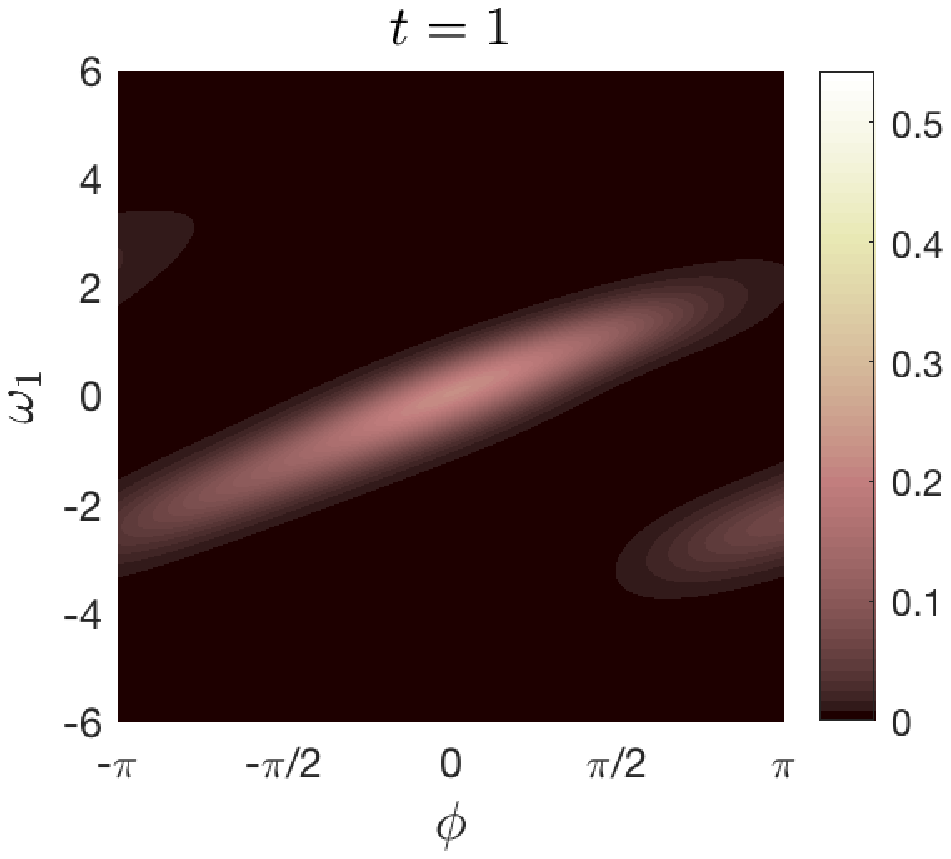}
\includegraphics[width = 0.328 \textwidth]{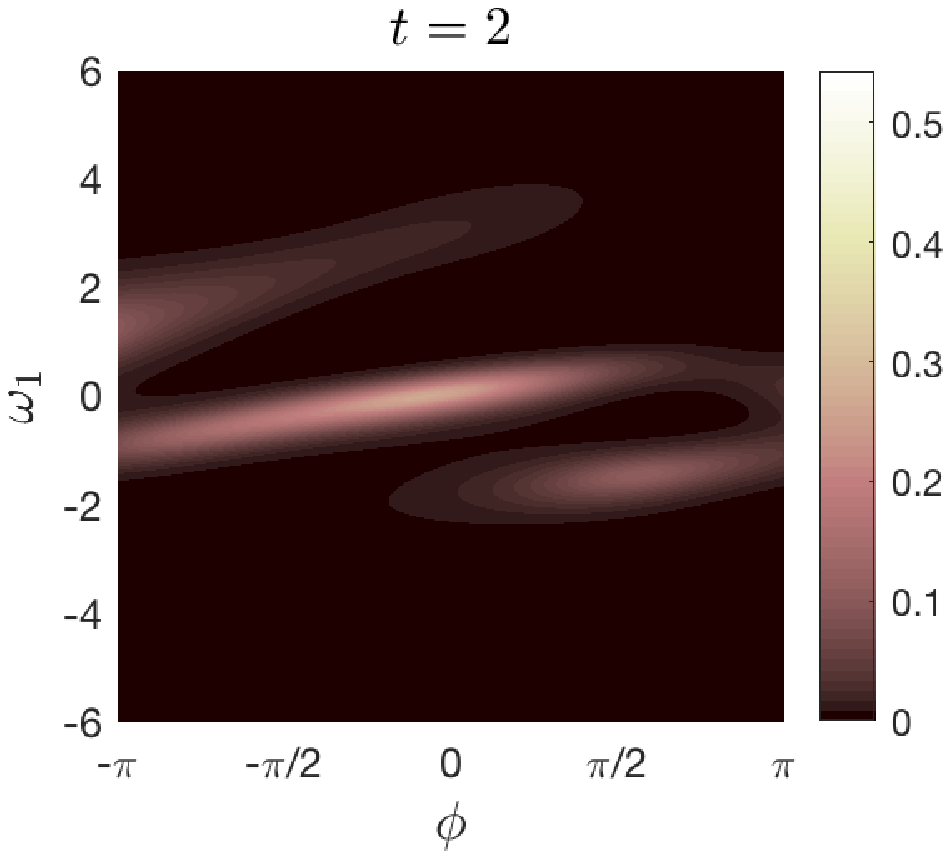}
\caption{$\rho_{\bm \theta} \left( \phi, \omega_1 \left| \theta = \psi = \omega_2 = \omega_3 = 0, 0 \leq \beta \leq 1/2 \right. ; t \right)$}
\label{fig: satellite pdf 1/2 evolution}
\end{subfigure}
\caption{Predicted conditional PDFs of the rigid body system with LQR control.}
\label{fig: satellite pdf evolution}
\end{figure*}

For this problem we train a NN with four hidden layers, each with 128 neurons. The weights in the regression loss term \cref{eq: data loss} are set to $w_i = 1$, $i = 1, \dots, | \mathcal D |$. In this problem the magnitude of the density grows almost exponentially in time, so we implement the transfer learning strategy from \cref{sec: transfer learning} with $N_t = 4$ time horizons,
$$
t_1 = 0.5,
\qquad
t_2 = 1,
\qquad
t_3 = 1.5,
\qquad
t_4 = 2,
$$
with corresponding PDE penalty weights,
$$
\lambda_1 = 1,
\qquad
\lambda_2 = 10,
\qquad
\lambda_3 = 10^5,
\qquad
\lambda_4 = 10^6.
$$

In the first three time horizons, we train for a single round on a fixed data set, $\mathcal D_1$, constructed from 500 sample trajectories. The collocation set includes $\mathcal D_1$ and an equal number of randomly sampled points, re-sampling every time horizon. For the final time horizon, $t_4$, we implement the adaptive learning algorithm with $\epsilon_{\rho} = 10^{-3}$, and $\epsilon_{\mathcal L} = 10^{-5}$. With these parameters, the convergence criteria (\ref{eq: norm test D}-\ref{eq: norm test C}) are satisfied after three rounds, observing data from a total of 984 trajectories evaluated at $K = 81$ time snapshots each.

In \cref{fig: satellite error} we plot the NRMSE \cref{eq: validation} of the NN at each time horizon, evaluated on a validation data set of 500 trajectories with $K = 101$ time snapshots each. This is compared to a direct solution attempt over the full time horizon, which turns out to be quite inaccurate. Thus we see that the transfer learning approach makes the proposed framework viable for approximating the rapidly growing density.

\Cref{fig: satellite pdf evolution} shows predicted conditional densities of $(\phi, \omega_1)$ for different values of $\beta$. In \cref{fig: satellite pdf 1 evolution}, we plot the case where $\beta = \bar \beta = 1$, the nominal value for which the control is designed. We can see that the density -- which is initially spread out in $\omega_1$ but tight in $\phi$ -- rotates so that $\omega_1$ converges to zero while $\phi$ spreads out. Then in \cref{fig: satellite pdf 1/2 evolution} we plot the case where $0 \leq \beta \leq 1/2$. We can clearly see that the density increases more slowly than when $\beta = 1$, indicating poorer stability. Still, overall the results are qualitatively similar, showing that the LQR controller has a large domain of attraction in $\omega_1$ and is fairly robust to changes in $\beta$, but does a poor job of stabilizing $\phi$. Similar results are easily obtained for $(\theta, \omega_2)$ and $(\psi, \omega_3)$ with the same NN. Finally, we observe that the NN predicts seemingly erroneous mass on the sides of the domain. But this is actually reasonable since $\phi$ and $\psi$ are periodic, and a cursory examination of the data reveals that many trajectories do in fact wrap around the periodic boundaries.

\section{Concluding remarks}

In this paper we have introduced a new data-driven approach to uncertainty propagation in nonlinear dynamic systems. Leveraging the approximating power of NNs, we model the probability density at any location $\bm x$ and time $t$, without discretizing the space-time domain. Numerical results suggest that the proposed method is competitive with other techniques for solving Liouville equations in moderate to high-dimensional systems. Different from most other approaches, our method may be more effective for approximating densities which are spatially complex, but with magnitude that doesn't vary too drastically over time. Thus we envision our method not as a replacement for existing techniques, but rather as filling a certain niche in uncertainty quantification: modeling densities which are \emph{not} characterized well by statistics like means and variances. Finally, we believe that this work represents an important step in developing physics-driven machine learning methods for optimal control in the space of probability distributions.

\bibliographystyle{siamplain}
\bibliography{bibliography}

\end{document}